\newif\ifarxiv
\arxivtrue

\ifarxiv
\documentclass[runningheads]{llncs}
\else 
\documentclass{llncs}
\fi

\usepackage[T1]{fontenc}
\usepackage{xspace}
\usepackage{todonotes}
\usepackage{graphicx}
\usepackage{amssymb}
\usepackage{enumitem}

\title{The thickness of fan-planar graphs is at most three%
\thanks{Research of Schlipf was supported by the Ministry of Science, Research and the Arts Baden-Württemberg (Germany).}}

\author{Otfried Cheong\inst{1}\orcidID{0000-0003-4467-7075}\and 
Maximilian Pfister\inst{2}\orcidID{0000-0002-7203-0669 }\and 
Lena Schlipf\inst{2}\orcidID{0000-0001-7043-1867} }

\institute{Institut f{\"u}r Informatik, Universit{\"a}t Bayreuth, \email{otfried.cheong@uni-bayreuth.de}
\and 
Wilhelm-Schickard-Institut f{\"u}r Informatik, Universit{\"a}t T{\"u}bingen, \email{maximilian.pfister@uni-tuebingen.de}, \email{lena.schlipf@uni-tuebingen.de}}

\usepackage{hyperref}
\hypersetup{%
  breaklinks,%
  colorlinks=true,%
  linkcolor=[rgb]{0.45,0.0,0.0},%
  citecolor=[rgb]{0,0,0.45}
}

\newcommand{\G}{\mathcal{G}}
\newcommand{\IG}{\mathcal{I}}
\newcommand{\Lo}{\mathcal{L}}
\newcommand{\R}{\mathcal{R}}
\newcommand{\Sp}{\mathcal{S}}
\newcommand{\base}{\hat{e}}

\newcommand{\refsimple}{(S)\xspace}
\newcommand{\refminimal}{(M)\xspace}
\newcommand{\reffani}{(I)\xspace}
\newcommand{\reffanii}{(II)\xspace}
\newcommand{\reffaniii}{(III)\xspace}

\begin{document}
\maketitle
\begin{abstract}
We prove that in any strongly fan-planar drawing of a graph~$\G$ the edges can be colored with
at most three colors, such that no two edges of the same color cross.  
This implies that the thickness of strongly fan-planar graphs is at most three.
If~$\G$ is bipartite, then two colors suffice to color the edges in this way.
\end{abstract}
\section{Introduction}
In order to visualize non-planar graphs, the research field of \emph{graph drawing beyond planarity} emerged as a
generalization of drawing planar graphs by allowing certain edge-crossing patterns.
Among the most popular classes of beyond-planar graph drawings are $k$-planar drawings, 
where every edge can have at most $k$~crossings, 
$k$-quasiplanar drawings, which do not contain $k$ mutually crossing edges, 
$RAC$-drawings, where the edges are straight-line segments that can only cross at right angles, 
fan-crossing-free drawings, where an edge is crossed by independent edges, and fan-planar drawings, 
where no edge is crossed by independent edges (see below for the full definition).
Many more details and further classes can be found in the recent
survey~\cite{10.1145/3301281}  and recent
book~\cite{Hong-Tokuyama} about beyond-planarity. 

\textbf{Thickness.}
A different notion of non-planarity was introduced by Tutte in~1963~\cite{tutte1963thickness}:
the \emph{thickness} of a graph is the minimum number of planar graphs into which the edges
of the graph can be partitioned. 
Subsequent research derived tight bounds on the thickness of complete and complete bipartite 
graphs~\cite{Alekseev_1976,beineke_harary_moon_1964}. 
For general graphs, however, it turns out to be hard to determine their thickness---even 
deciding if the thickness of a given graph is two is already NP-hard~\cite{mansfield_1983}. 
More details about thickness and related concepts can be found in~\cite{thickness-survey}.

The first connection between thickness and beyond-planarity was laid in~1973, 
when Kainen~\cite{kainen1973} studied the relationship between the thickness of a graph and
$k$-planarity. 
He observed that a $k$-planar graph~$\G$ has thickness at most~$k+1$. 
This follows from the fact that~$\G$ admits by definition a drawing
where every edge crosses at most $k$~other edges, hence the edge intersection graph 
(refer to Section~\ref{sec:def} for a formal definition) has maximum degree $k$
and admits a vertex coloring using at most $k+1$~colors.
For other beyond-planar graph classes, however, the thickness is not so closely related
to a coloring of the edge intersection graph.
For $3$-quasiplanar graphs, for instance, there exist $3$-quasi-planar drawings 
whose edge intersection graph cannot be colored with a bounded number of colors~\cite{pawlik2014triangle}, 
yet the thickness of $3$-quasi-planar graphs is at most~seven. 

The most general tool to determine the thickness of a graph~$\G$ is to compute its \emph{arboricity}, 
which is the minimum number of forests into which the edges of~$\G$ can be partitioned. 
Since a forest is planar, the arboricity is an upper bound on the thickness, 
and since any planar graph has arboricity at most three, 
the arboricity of~$\G$ is at most three times its thickness. 
By the Nash-Williams theorem~\cite{10.1112/jlms/s1-36.1.445}, the arboricity of~$\G$ is~$a(\G) = \max_S{\frac{m_S}{n_S-1}}$,
where~$S$ ranges over all subgraphs of~$\G$ with~$m_S$ edges and~$n_S$ vertices. 
For typical beyond-planar graph classes, a subgraph of a graph is contained inside the same class,
and so a linear bound on the number of edges of graphs in the class implies a constant bound on the thickness.
This insight actually improves the result of Kainen for large enough~$k$, 
as $k$-planar graphs have at most~$3.81\sqrt{k}n$ edges~\cite{ackerman2015topological},
and therefore the thickness of $k$-planar graphs is~$\mathcal{O}(\sqrt{k})$.

\textbf{Fan-planar graphs.}
We consider the class of \emph{fan-planar} graphs.
This class was introduced in 2014 by Kaufmann and Ueckerdt~\cite{DBLP:journals/corr/KaufmannU14} as graphs that admit a \emph{fan-planar drawing}, 
which they define by the requirement that for each edge~$e$, the edges crossing~$e$ have a common 
endpoint on the same side of~$e$.  
This can be formulated equivalently by two forbidden patterns, patterns~(I) and~(II) in Fig.~\ref{fig:fanplanar}: 
one is the configuration where~$e$ is crossed
by two independent edges, the other where~$e$ is crossed by incident edges with the common endpoint
on different sides of~$e$.
Fan-planarity has evolved into a popular subject of study with many related publications, 
a good overview of which can be found in the recent survey article by Bekos and Grilli~\cite{Bekos2020}.

And now it get's complicated: the recent journal version~\cite{DBLP:journals/combinatorics/0001U22} 
of the 2014 paper~\cite{DBLP:journals/corr/KaufmannU14} gives a more restricted definition of fan-planarity, 
where also pattern~(III) of Fig.~\ref{fig:fanplanar} is forbidden.  
The restriction is necessary to allow the proof for 
the bound~$5n-10$ on the number of edges of a fan-planar graph to go through.  
Patterns~(II) and~(III) can formally be defined as follows: 
if one removes the edge~$e$ and the two edges intersecting~$e$ from the plane, 
we are left with two connected regions called cells, a bounded cell and an unbounded cell. 
In pattern~(II), one endpoint of~$e$ lies in the bounded cell, in pattern~(III), 
both endpoints lie in the bounded cell.
Both patterns are forbidden, so both endpoints of~$e$ must lie in the unbounded cell.

To distinguish this definition from the ``classic'' definition, 
let's call a fan-planar drawing without pattern~(III) a \emph{strongly fan-planar} drawing,
and a graph with such a drawing a \emph{strongly fan-planar} graph.
Note that---quite unusual for a topological drawing style---strongly fan-planar drawings are
only defined in the plane, the definition does not work for drawings on the sphere!
Clearly, strongly fan-planar graphs are also fan-planar (in the ``classic'' sense), 
so most of the literature applies to this restricted class as well.

\textbf{Our Contribution.}
We establish the first result on the thickness of a beyond-planar graph class 
that is stronger than bounds based only on the density or arboricity of the class.
We show that the edges of a simple, strongly fan-planar drawing of a graph~$\G$ 
can be colored with three colors, 
such that no two edges of the same color cross, implying that~$\G$ has thickness at most three.
We also show that if~$\G$ is bipartite, then two colors suffice for the same coloring, 
and therefore bipartite strongly fan-planar graphs have thickness two (and this bound is tight unless~$\G$ is planar).
For comparison, the upper bound implied by density (and therefore arboricity) of strongly fan-planar graphs is 
five~\cite{DBLP:journals/combinatorics/0001U22} 
(and four for bipartite strongly fan-planar graphs~\cite{angelini_et_al:LIPIcs:2018:9976}).
Our proof relies on a complete characterization of chordless cycles in the intersection graph,
which should be of independent interest for the study of fan-planar drawings.


\section{Preliminaries}
\label{sec:def}

Throughout the paper, we assume that any graph~$\G$ and its corresponding drawing is \emph{simple}, 
that is, $\G$ has no self-loops or multiple edges, adjacent edges cannot cross, 
any two edges are allowed to cross at most once, 
and crossing points of distinct pairs of edges do not coincide.
We will refer to all of these as property~\refsimple.
(Non-simple fan-planar graphs where discussed by Klemz et al.~\cite{Klemz-Knorr-Reddy-Schroeder}.)
\begin{figure}[t]
  \centerline{\includegraphics[page=8,scale=1]{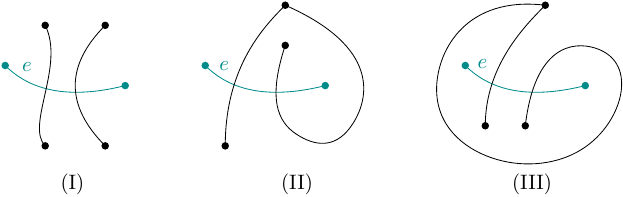}}
  \caption{The three forbidden configurations in strongly fan-planar drawings.}
  \label{fig:fanplanar}
\end{figure}
The \emph{intersection graph} $\IG$ of a drawing $\Gamma$ of graph~$\G$ has a vertex for 
every edge of~$\G$, and two vertices are connected in~$\IG$ if the corresponding edges cross in $\Gamma$.

Let $\G$ be a strongly fan-planar graph  with a fixed, strongly fan-planar drawing~$\Gamma$.
(Throughout the proofs, we will simply say ``fan-planar'' with the understanding that we require 
strong fan-planarity.)
This means that none of the forbidden configurations shown in Fig.~\ref{fig:fanplanar} occur.
For ease of reference, we will refer to the three forbidden patterns as 
properties~\reffani, \reffanii, and~\reffaniii.
If an edge~$e$ is crossed by more than one edge, all edges
crossing~$e$ share a common endpoint---we call this point the~$\emph{anchor}$ of~$e$.
Throughout the paper, we will use the letter~$v$ to denote the anchor of the edge denoted by~$e$, 
so that~$v_i$ is always the anchor of~$e_i$.

We will discuss \emph{chordless} cycles in the intersection graph~$\IG$ of~$\G$, 
that is, cycles without diagonals. 
In~$\G$, a chordless cycle $C$~corresponds to a sequence of 
edges~$e_1, \dots, e_k$, such that $e_i$ and~$e_{i+1}$ intersect, 
but there are no other intersections between the edges of $C$---we will refer to this
as property~\refminimal. (Throughout the paper, 
all arithmetic on indices is modulo the size of the cycle.)
Fixing a chordless cycle~$C$, we define~$x_i$ as the intersection point 
between~$e_{i-1}$ and~$e_i$. We let $a_i$ and $b_i$ be the endpoints
of~$e_i$, such that $a_i x_i x_{i+1} b_i$ appear in this order on~$e_i$.
We will call~$a_i$ the \emph{source},~$b_i$ the \emph{target} of~$e_i$---but keep in mind that
this orientation is only defined with respect to~$C$.
Let $\base_i$ be the oriented segment of~$e_i$ from~$x_i$
to~$x_{i+1}$---we will call~$\base_i$ the \emph{base} of~$e_i$.
If we concatenate the bases~$\base_i$ in order, we
obtain a closed loop that we call~$\Lo$.  Since by~\refminimal two
bases~$\base_i$, $\base_j$ do not intersect, 
$\Lo$ is a Jordan curve and partitions the plane into two regions.
Since $\G$ is fan-planar, edges~$e_i$ and~$e_{i+2}$ share an endpoint, 
namely the anchor~$v_{i+1}$ of~$e_{i+1}$.
We will use~$\G_C$ to denote the subgraph of~$\G$ consisting of only the edges of~$C$, 
with the same embedding.

\begin{figure}[h!]
  \centerline{\includegraphics[page=3]{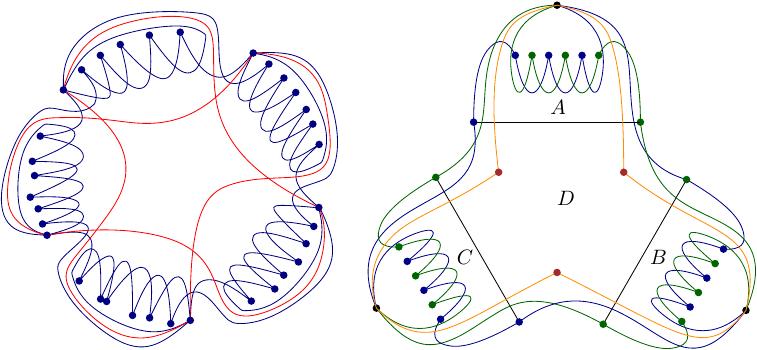}}
  \caption{Chordless cycles can intersect and share edges.}
  \label{fig:cyclesintersect}
\end{figure}

\section{Characterizing chordless cycles}
\label{sec:cycle}
In this section we characterize chordless cycles in~$\IG$.  
In the next section we then study how chordless cycles can interact, 
and we will be able to break all odd cycles simultaneously by coloring a carefully chosen set of edges
with one color.
Fig.~\ref{fig:cyclesintersect} shows two examples of strongly fan-planar graphs 
whose intersection graphs have several chordless cycles 
that cross and share edges.  The graph on the left has 32~distinct chordless cycles, 
as for each of the five red edges we can instead traverse the blue edges ``behind.''
In the graph on the right, the boundaries of the faces labeled~$A, B, C$ are loops of chordless cycles of length~11.
There is a chordless cycle of length~9 surrounding face~$D$, there is a chordless cycle of length~30 surrounding all four faces,
and there are three chordless cycles of length~23 surrounding~$ABD$,~$BCD$, and~$CAD$, respectively.

\begin{figure}[tb]
  \centerline{\includegraphics[page=2,width=\textwidth]{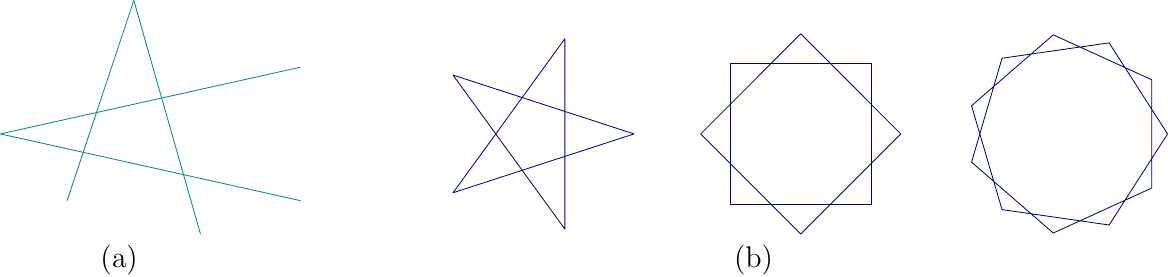}}
  \caption{(a) The only possible cycle of length four, and (b) fully canonical cycles for~$k=5, 8, 11$.}
  \label{fig:regular}
\end{figure}

It is easy to see using~\reffani and~\refsimple that~$\IG$ cannot have cycles of length three, 
while cycles of length four have a unique shape, see Fig.~\ref{fig:regular}(a).
It remains to study the structure of chordless cycles of length at least five.

We will call an edge~$e_i$ of a chordless cycle \emph{canonical} if the anchor~$v_i$ 
of~$e_i$ is the target of $e_{i-1}$ and the source of~$e_{i+1}$, that is, if~$b_{i-1} = v_i = a_{i+1}$. 
If, in addition, no other edge of~$C$ is incident to~$v_i$---that is, 
if $v_i$ has degree two in~$\G_C$---then~$e_i$ is \emph{strictly canonical}.

\begin{figure}[htb]
  \centerline{\includegraphics[page=1]{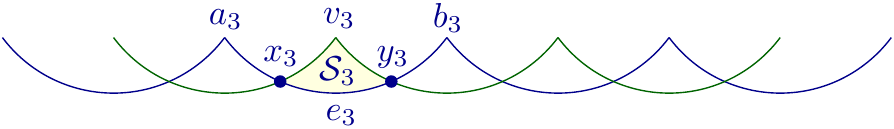}}
  \caption{A sequence of canonical edges.}
  \label{fig:canonical}
\end{figure}

\begin{figure}[h]
  \centerline{\includegraphics[page=5]{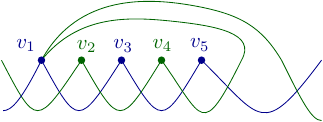}}
  \caption{The anchors of canonical edges can coincide, but only after at least four distinct anchors.}
  \label{fig:strictcanonical}
\end{figure}

Fig.~\ref{fig:canonical} shows a sequence of canonical edges. 
Note that some of the endpoints of the edges of this 
sequence can coincide, see Fig.~\ref{fig:strictcanonical}. 
For a canonical edge~$e_i$, we will call the ``triangle''~$\Sp_i$ with corners $x_i$, 
$x_{i+1}$, $v_i$ and
bounded by the edges~$e_{i-1}$, $e_{i}$, and~$e_{i+1}$ the \emph{spike} of~$e_i$, 
see the shaded region in Fig.~\ref{fig:canonical}.

\begin{lemma}
  \label{lem:strict}
  Let $e_1,e_2,e_3,e_4,e_5$ be five consecutive canonical edges of a chordless cycle.  Then the anchors of the five edges are distinct. 
\end{lemma}
\begin{proof}
Since $v_2$ is an endpoint of~$e_1$, $v_2 \neq v_1$ holds. By construction, $v_1$ and $v_3$ are the endpoints of~$e_2$, so $v_3 \neq v_1$.  If $v_1 = v_4$, then~$e_2$ and~$e_3$ share an endpoint and intersect, a contradiction to~\refsimple. Finally, if $v_1 = v_5$, then both~$e_2$ and~$e_4$ connect~$v_1$ and~$v_5$ and are therefore identical. Using analogous arguments one establishes the remaining inequalities. \qed
\end{proof}

\begin{figure}[h]
  \centerline{\includegraphics[page=14]{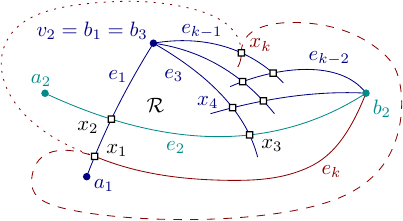}}
  \caption{The cycle cannot be closed.}
  \label{fig:bisequal}
\end{figure}

\begin{lemma}\label{lem:bnotequal}
  Let $C = (e_1,\dots e_k)$ be a chordless cycle in~$\IG$ with~$k \geq 5$. Then, $b_{i-1} \neq b_{i+1}$ and
  $a_{i-1} \neq a_{i+1}$ for~$1 \leq i \leq k$.
\end{lemma}
\begin{proof}
Suppose for a contradiction that there is an index~$i$ such that $b_{i-1} = b_{i+1}$. By renumbering we can assume that~$b_1 = v_2 = b_3$.
This implies that $x_4$ lies on the segment $x_3v_2$, and $x_1$ lies on the segment~$a_1x_2$, see Fig.~\ref{fig:bisequal}.
Consider~$e_4$, which contains~$x_4$ and is incident to~$v_3$. 
Since by~\refminimal $e_4$ does not cross~$e_1$, we have $v_3 \neq a_2$ as otherwise
$e_3$ forms forbidden configuration~\reffanii with~$e_2$ and~$e_4$.
Consider the region $\R$ bounded by $v_2 x_2 x_3$. 
By~\refminimal the two boundary segments $v_2x_2$ and $x_2x_3$ do not cross any edge of~$C$, 
and the segment $x_3v_2$ only crosses~$e_4$.
It follows that~$e_4$ is incident to~$b_2$ and its other endpoint lies in~$\R$.
If~$x_5$ lies inside~$\R$, then the closed curve~$\Lo$ has to intersect the boundary of~$\R$ to return to~$x_1$, a contradiction.
So~$x_5$ lies on~$x_4 b_2$, and therefore~$b_4 = b_2$. 
The same argument now implies that~$b_5 = b_3 = b_1$, $b_6 = b_4 = b_2$, and so on.
The final edge~$e_k$ 
must contain~$x_1$, so it cannot be incident to~$b_1$ by~\refsimple, and so~$k$ is even, $e_k$ is incident to~$b_2$, and~$x_k$ lies on the edge~$e_{k-1}$ incident to~$b_1$.
Since $x_1$ lies on~$e_k$ between~$x_k$ and~$b_k = b_2$, there are two possibilities for drawing~$e_k$, 
shown dashed and dotted in Fig.~\ref{fig:bisequal}.
The dotted version violates~\reffaniii for~$e_{k-1}$ (with~$e_{k-2}$ and~$e_k$),
the dashed version violates~\reffaniii for~$e_k$ (with~$e_1$ and~$e_{k-1}$), a contradiction.

Assume now that~$a_{i-1} = a_{i+1}$ for some~$i$, and let~$C' = (e_k,\dots,e_1)$ be the chordless cycle 
obtained by reversing~$C$. Reversing the direction flips~$a_i$ and~$b_i$ for each edge, so~$C'$ now has 
an index~$j$ such that~$b_{j-1} = b_{j+1}$, and we already showed this cannot be the case. \qed
\end{proof}

\begin{corollary}
  \label{cor:noncanon}
  Let $e_i$ be a non-canonical edge of a chordless cycle~$C$. 
  Then $a_{i-1} = v_i = b_{i+1}$.
\end{corollary}
\begin{proof}
  The anchor~$v_i$ is a common endpoint of~$e_{i-1}$ and~$e_{i+1}$.  This cannot be~$b_{i-1}$, 
  because $e_i$ is not canonical and~$b_{i-1} \neq b_{i+1}$ by Lemma~\ref{lem:bnotequal}. 
  So~$v_i = a_{i-1}$, and since~$a_{i-1} \neq a_{i+1}$ by Lemma~\ref{lem:bnotequal} again, 
  we have~$v_i = a_{i-1} = b_{i+1}$. \qed
\end{proof}

We say that the chordless cycle~$C$ is \emph{fully canonical} if all its edges are canonical.
Such cycles can be created for any~$k \geq 5$ by taking the corners of a regular~$k$-gon and
connecting every other corner with an edge, see
Fig.~\ref{fig:regular}(b).
Note that such a fully canonical cycle corresponds to a single closed trail of length~$k$
in~$\G$ for odd~$k$, but to two closed trails of length~$k/2$ for even~$k$.  
These closed trails in~$\G$ are not necessarily cycles, 
as the anchors of the edges of a fully canonical cycle can coincide, 
see e.g.~Fig.~\ref{fig:cyclesintersect}(left).

\begin{lemma}
  \label{lem:canonicaledges}
  A chordless cycle of length at least five has at least four consecutive canonical edges.
\end{lemma}
\begin{proof}
Let $C = (e_1,\dots e_k)$ be a chordless cycle of length $k\geq 5$ such that~$e_1,\dots,e_m$ is a longest consecutive sequence
of canonical edges in~$C$, and assume~$m < 4$.  We distinguish four different cases according to the value of~$m$. 
\begin{figure}[h]
  \centerline{\includegraphics[page=13]{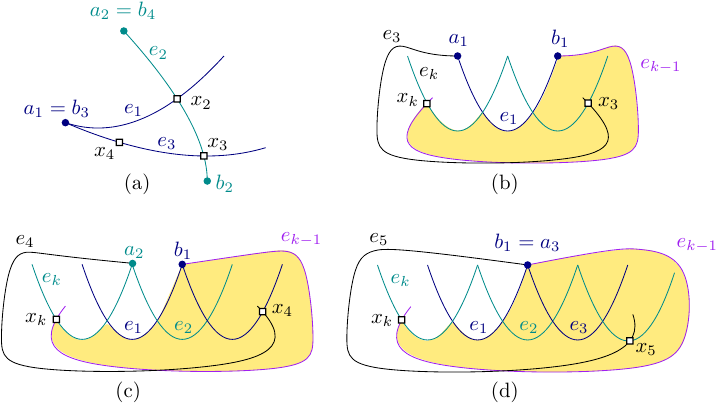}}
  \caption{(a) $m=0$, (b) $m=1$, (c) $m=2$ and (d) $m=3$.}
  \label{fig:fourcanonical}
\end{figure}

If~$m = 0$, that is, $C$ does not have any canonical edge:
By Corollary~\ref{cor:noncanon}, $a_1 = v_2 = b_3$ and~$a_2 = v_3 = b_4$. 
This implies that~$x_4$ lies on the segment~$b_3x_3$, and~$e_4$ must connect~$b_4 = a_2$ with~$x_4$
without intersecting either~$e_1$ or~$e_2$, creating forbidden configuration~\reffanii,
see Fig.~\ref{fig:fourcanonical}(a).

If~$m = 1$: 
Since~$e_2$ and~$e_k$ are not canonical, by Corollary~\ref{cor:noncanon} $a_1 = b_3$ and $b_{1}=a_{k-1}$.
By~\refminimal, $e_3$ does not cross~$e_k$, and~$e_2$ does not cross~$e_{k-1}$. 
Therefore,~$x_k \in e_{k-1}$  and~$x_3 \in e_3$ implies
that we have the situation of Fig.~\ref{fig:fourcanonical}(b) (by \reffanii and \reffaniii).  But then~$e_3$ and~$e_{k-1}$ cross, so we must have~$k=5$. 
Since~$e_4$ crosses~$e_3$ and~$e_5$, $e_3$ and $e_5$ must have a common endpoint.
By~\refsimple, $e_5$~cannot be incident to~$a_1$, and~$e_3$ cannot reach either endpoint of~$e_5$ without introducing a crossing
that contradicts~\refsimple.

If~$m = 2$: 
Since~$e_3$ and~$e_k$ are not canonical, by Corollary~\ref{cor:noncanon} $a_2 = b_4$ and $b_{1}=a_{k-1}$.
Since~$e_4$ and~$e_k$ have a common endpoint, they cannot cross, so~$k > 5$.  Also~$e_{k-1}$ and~$e_3$ do not cross.
Therefore~$x_k \in e_{k-1}$ and~$x_4 \in e_4$ implies
that we have the situation of Fig.~\ref{fig:fourcanonical}(c). So~$e_4$ and~$e_{k-1}$ cross and we have~$k = 6$.
However, $e_4,e_5,e_6$ (and also~$e_3,e_4,e_5$) form forbidden configuration~\reffanii.

If~$m = 3$: 
Since~$e_4$ and~$e_k$ are not canonical, by Corollary~\ref{cor:noncanon} $a_3 = b_5$ and $b_{1}=a_{k-1}$, so, since~$e_2$ is canonical,
we have~$a_{k-1} = b_1 = a_3 = b_ 5$.
Since~$e_5$ and~$e_{k-1}$ share this endpoint, they cannot cross, so~$k \neq 5$.
If~$k = 6$, then~$a_5 = b_5$, a contradiction, so we have~$k \geq 7$.
It follows that~$e_5$ does not cross~$\{e_1, e_2, e_3, e_k\}$,
and~$e_{k-1}$ does not cross~$\{e_1, e_2, e_3, e_4\}$.
Therefore~$x_k \in e_{k-1}$ and~$x_5 \in e_5$ implies
that we have the situation of Fig.~\ref{fig:fourcanonical}(d). 
But this requires~$e_5$ and~$e_{k-1}$ to cross, a contradiction to the above. \qed
\end{proof}

\begin{theorem}
  \label{thm:structure}
  If a chordless cycle of length~$k \geq 5$ is not fully canonical, 
  then $k \geq 9$, edges~$e_1,\dots,e_{k-1}$ are canonical, 
  anchors~$v_2 = v_{k-2}$ coincide so that
  $b_1 = a_3 = b_{k-3} = a_{k-1}$, and $b_{k-1}$ and~$a_1$ 
  are vertices of degree one in~$\G_C$.
\end{theorem}
\begin{proof}
Consider a chordless cycle~$C = (e_1,\dots,e_k)$ of length $k \geq 5$ that is not fully canonical, and
such that~$e_1,\dots,e_{m-1}$ is a longest sub-sequence of~$C$ consisting of canonical edges. 
By Lemma~\ref{lem:canonicaledges} we have~$m \geq 5$. Since~$e_m$ is not canonical, we have~$a_{m-1} = b_{m+1}$ by Corollary~\ref{cor:noncanon}.
By definition,~$e_{m+1}$ crosses~$e_m$ in~$x_{m+1}$, and by~\refminimal it crosses no other edge~$e_{i}$ with $i \neq m+2$.
Let~$\R$ be the region enclosed by the base~$\base_m$ and the edges~$e_{m-1}$ and~$e_{m+1}$. By~\reffanii and~\reffaniii,
$\R$ contains no endpoint of~$e_m$, and so we have the situation of Fig.~\ref{fig:endcanonical}(a).
\begin{figure}[h]
  \centerline{\includegraphics[page=4,width=\textwidth]{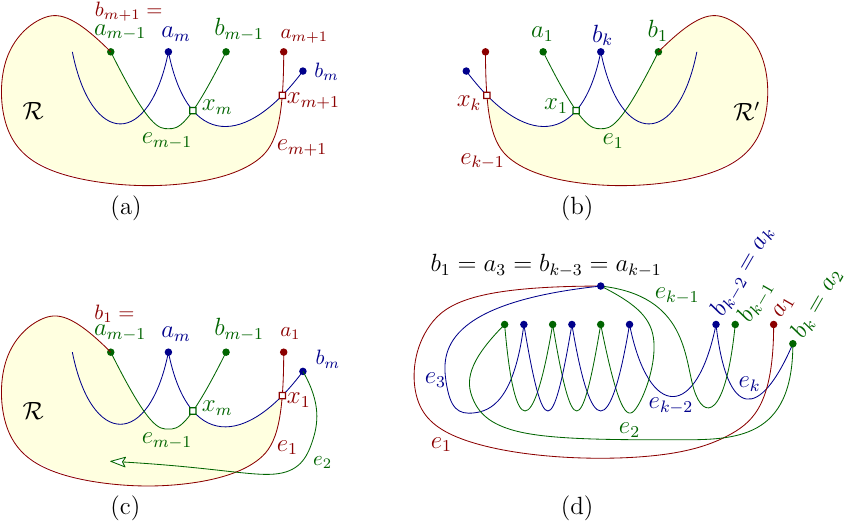}}
  \caption{Proof of Theorem~\ref{thm:structure}.}
  \label{fig:endcanonical}
\end{figure}
The loop~$\Lo$ of~$C$ lies entirely in~$\R$, and only the bases~$\base_{m-1}$, $\base_m$, $\base_{m+1}$ lie on the boundary of~$\R$.

We want to show that~$k = m$, so we assume for a contradiction that~$k > m$. 
If~$k = m+1$, then~$e_k = e_{m+1}$ in Fig.~\ref{fig:endcanonical}(a).
Since~$e_k$ is not canonical, we then have~$a_m = b_1$, 
and~$e_1$ intersects~$e_k = e_{m+1}$ between~$x_{m+1}$ and~$b_{m+1} = a_{m-1}$.
Since~$e_1$ cannot cross~$e_{m-1}$ by~\refminimal, that violates either~\reffanii or~\reffaniii.

If~$k = m+2$, then~$e_{k-1} = e_{m+1}$ in Fig.~\ref{fig:endcanonical}(a),
$x_k$~lies on~$e_{m+1}$ on the boundary of~$\R$, 
and~$\base_k$ lies (except for its endpoint~$x_k$) in the interior of~$\R$. 
But since~$e_k$ is not canonical, we have~$a_{m+1} = b_1$,
and since~$e_1$ does not cross the boundary of~$\R$, it cannot contain~$x_1$ in the interior of~$\R$.

We assume next that~$k > m+2$, which implies that~$\base_k$ lies in the interior of~$R$.
Since~$e_k$ is not canonical, we have~$a_{k-1} = b_1$. 
Symmetrically to the argument above, 
the edge~$e_{k-1}$ must be such that the region~$\R'$ formed by~$\base_{k}$, $e_{k-1}$, 
and~$e_1$ contains no endpoint of~$e_k$, so we are in the situation of Fig.~\ref{fig:endcanonical}(b).
Again, the loop~$\Lo$ lies in~$\R'$, with only~$\base_{k-1}, \base_k, \base_1$ on the boundary of~$\R'$.
In particular,~$\base_m$ lies in the interior of~$\R'$.

Since~$\base_m$ is in the interior of~$\R'$ but on the boundary of~$\R$, 
while~$\base_k$ is in the interior of~$\R$ but on the boundary of~$\R'$, 
the two regions cannot be nested, and their boundaries must intersect.
The boundary of~$\R$ consists of~$e_{m-1}$, $\base_m$, and~$e_{m+1}$,
the boundary of~$\R'$ consists of~$e_{k-1}$, $\base_k$, and~$e_1$, 
so by~\refminimal the only possible edge crossing occurs when~$k = m+3$, 
so that~$e_{m+1}$ and~$e_{k-1}$ can cross. 
Since two intersecting closed curves must intersect an even number of times,
we must in addition have that the vertices~$a_{m-1}$ of~$\R$ and~$b_1$ of~$\R'$ coincide,
but then~$e_{m+1}$ and~$e_{k-1}$ have a common endpoint and cannot cross at all by~\refsimple.

It follows that our assumption that~$k > m$ is false, and so~$k = m$.  
Relabeling the edges in Fig.~\ref{fig:endcanonical}(a) we obtain Fig.~\ref{fig:endcanonical}(c).
Since~$e_1$ is canonical, we have~$b_m = a_2$, so~$e_2$ starts in~$b_m$ and enters~$\R$ through~$e_1$,
see Fig.~\ref{fig:endcanonical}(c).
Since~$e_2$ cannot cross either~$e_{m-1}$ or~$e_m$, its other endpoint~$b_2$ lies in the interior of~$\R$.
Now it remains to observe that since~$e_2, \dots, e_{k-2}$ are canonical, 
we end up with the situation shown in Fig.~\ref{fig:endcanonical}(d).

Since~$v_2 = b_1 = a_{k-1} = v_{k-2}$, Lemma~\ref{lem:strict} implies that $k-2 \geq 7$, so~$k \geq 9$. \qed
\end{proof}

\begin{figure}[h]
  \centerline{\includegraphics[page=3]{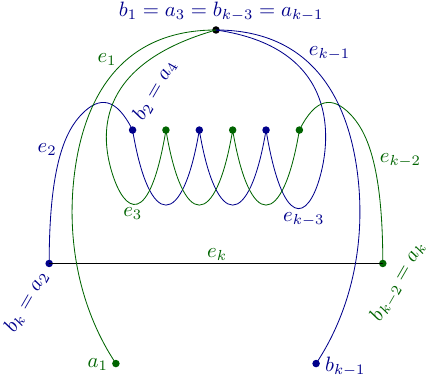}}
  \caption{A non-canonical chordless cycle.}
  \label{fig:noncanonical}
\end{figure}
In Fig.~\ref{fig:noncanonical} we draw again the non-canonical chordless cycle 
of~Fig.~\ref{fig:endcanonical}(d), showing the symmetry in the characterization.
The reader may enjoy determining which of the chordless cycles 
in~Fig.~\ref{fig:cyclesintersect} are fully canonical, 
and which edges are the non-canonical ones of the other cycles.

It is easy to see that when~$k$ is odd, then the graph~$\G$ contains a closed trail of
length~$k-2$ consisting of all edges of~$C$ except for~$e_{k-1}$
and~$e_1$, namely the closed trail~$e_2 e_4 e_6 \dots e_{k-3} e_3 e_5 \dots
e_{k-2} e_k$.
In Fig.~\ref{fig:endcanonical}(d) we have~$k = 11$, so there is a closed trail of
length~$9$. Fig.~\ref{fig:noncanonical} shows the smallest possible non-canonical chordless cycle: 
here~$k = 9$, and so there is a closed trail of length~$7$.
We obtain the following corollary: 
\begin{corollary}
  The edges of a strongly fan-planar drawing of a bipartite graph~$\G$ can be colored using two colors such that no edges of the same color cross.
  As a consequence, a bipartite, strongly fan-planar graph has thickness at most two.
\end{corollary}
\begin{proof}
  We show that~$\IG$ is bipartite.
  Assume otherwise:  then~$\IG$ has an odd cycle, which contains a
  chordless odd cycle~$C$ of length~$k$.  If $C$ is fully canonical, then the edges of~$C$ form an odd cycle of length~$k$ in~$\G$, so~$\G$ is not bipartite.
  If $C$ is not fully canonical, then Theorem~\ref{thm:structure} implies
  that~$\G$ has an odd cycle of length~$k-2$, again a contradiction. \qed
\end{proof}

The bound on the thickness is tight: every bipartite, strongly fan-planar graph that is not planar 
has thickness exactly two, an example being~$K_{3,3}$.

In the remainder of this section, 
we study the regions induced by the loop ~$\Lo$ and the spikes~$\Sp_i$ in a little more detail, 
in particular how they can be crossed by \emph{other} edges, that is, edges not part of~$C$. 
The edges of~$C$ intersect only in the corners of the loop~$\Lo$.
The loop~$\Lo$ partitions the plane into two regions.  
If~$C$ is fully canonical, then one of these regions is empty in~$\G_C$---this could be either
the bounded or the unbounded region delimited by~$\Lo$.
If~$C$ has a non-canonical edge, 
then from Theorem~\ref{thm:structure} it follows that the \emph{bounded} region 
delimited by~$\Lo$ is empty in~$\G_C$ 
(here it plays a role that property~\reffaniii cannot be used for drawings on a sphere).
In both cases, we will denote the region delimited by~$\Lo$ that is empty in~$\G_C$ by~$\Lo$ as well.
Adjacent to~$\Lo$ is, for each base~$\base_i$ of a canonical edge~$e_i$,
the spike~$\Sp_i$, which is itself an empty region~in~$\G_C$. 

\begin{figure}[htb]
  \centerline{\includegraphics[page=6,width=\textwidth]{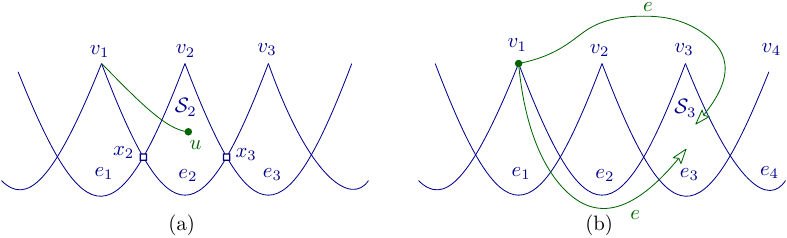}}
  \caption{(a) $u$ must lie in~$\Sp_2$ or $e = e_2$, and (b) $e$ starting in~$v_1$ cannot enter~$\Sp_3$.}
  \label{fig:vertex-in-spike}
\end{figure}
\begin{lemma}
  \label{lem:vertex-is-in-spike}
  Let $e_1,e_2,e_3$ be three consecutive canonical edges of a chordless cycle. Let  $e = (v_1, u)$ be an edge that crosses~$e_1$ in the relative interior of the segment~$x_2v_2$. Then $u$ is contained in~$\Sp_2$.
\end{lemma}
\begin{proof}
Traversing~$e$ from $v_1$ to $u$, $e$ enters $\Sp_2$ by crossing the segment $x_2 v_2$,
see Fig.~\ref{fig:vertex-in-spike}(a).
Since~$e$ cannot cross~$e_2$ and cannot cross~$e_1$ again by~\refsimple, 
$e$ could leave~$\Sp_2$ only through the segment~$v_2x_3$, by crossing~$e_3$. 
But then $v_3$ would have to be an endpoint of~$e$. 
By Lemma~\ref{lem:strict}, $v_3 \neq v_1$, so $u = v_3$, but then $e = e_2$, 
a contradiction since~$e_2$ does not intersect the relative interior of~$x_2 v_2$. \qed
\end{proof}

\begin{lemma}
  \label{lem:vertex-in-next-spike}
  Let $e_1,e_2,e_3,e_4$ be four consecutive canonical edges of a chordless cycle. 
  Let  $e = (v_1, u)$ be an edge incident to~$v_1$. Then~$e$ does not enter the interior of~$\Sp_3$.
\end{lemma}
\begin{proof}
By~\refsimple $e$~cannot cross~$e_2$, so to enter~$\Sp_3$, it would have to cross either~$e_3$ or~$e_4$, see Fig.~\ref{fig:vertex-in-spike}(b).
If~$e$ crosses~$e_3$, then it must be incident to~$v_3$, and since~$v_3 \neq v_1$, 
that means~$e = (v_1, v_3) = e_2$, which does not intersect the interior of~$\Sp_3$.
If~$e$ crosses~$e_4$, then it must be incident to~$v_4$, and since~$v_4 \neq v_1$,
that means~$e = (v_1, v_4)$. But~$v_4$ lies outside~$\Sp_3$, so~$e$ would have to cross~$e_4$ again
to reach~$v_4$, a contradiction to~\refsimple. \qed
\end{proof}

\begin{lemma}
 \label{lem:e-cross-l}
 Let~$C$ be a chordless cycle of length at least five, 
 and let~$e$ be an edge not part of~$C$ such that~$e$ intersects the loop~$\Lo$ of~$C$. 
 Then~$e$ starts in the anchor~$v_i$ of a canonical edge~$e_i$ of~$C$, 
 passes through the spike~$\Sp_i$, 
 crosses the base~$\base_i$, and either 
 (1)~ends in~$\Lo$; 
 or (2)~crosses the base~$\base_j$ of another canonical edge~$e_j$ of~$C$, 
 then passes through~$\Sp_j$ and terminates in~$v_j$; or 
 (3)~crosses the base of the non-canonical edge~$e_j$ of~$C$, never enters~$\Lo$ again, 
 and terminates in a vertex that is not a vertex of~$C$. 
 The second and third case can only happen if~$e_i$ and~$e_j$ share an endpoint. 
 In particular, (2) cannot happen when~$e_{i-1}$ and~$e_{i+1}$ are strictly canonical, 
 and (3) implies that~$i \in \{j-2, j+2\}$.
\end{lemma}
\ifarxiv
The proof can be found in Appendix~\ref{sec:e-cross-l}.
\fi

\newpage
\section{Coloring with three colors}
\label{sec:colorwiththree}

In this section we will show our main theorem:
\begin{theorem}
  \label{thm:main}
  In every strongly fan-planar drawing of a graph~$\G$ there is a set~$S$ of edges 
  such that (1)~$S$ is independent in~$\IG$, that is, no two edges in~$S$ cross; 
  and (2) every odd cycle in~$\IG$ contains an edge in~$S$.
\end{theorem}
The theorem immediately implies the following:
\begin{corollary}\label{cor:3colorable}
  The edges of a strongly fan-planar drawing of a graph~$\G$ can be colored using three colors such that
  no two edges of the same color cross.
  As a consequence, a strongly fan-planar graph has thickness at most three.
\end{corollary}
\begin{proof}
  Pick the set~$S$ of edges according to Theorem~\ref{thm:main} and color them with the first color.  Then~$\IG \setminus S$ contains no odd cycle and is therefore bipartite, and can be colored with the remaining two colors. \qed
\end{proof}
We construct the set~$S$ for the proof of Theorem~\ref{thm:main} using the following lemma.
\begin{lemma}\label{lem:ground}
    Let $C$ be a chordless cycle of length at least five. Then $C$ contains an edge~$e_i$
    such that 
    \begin{itemize}[nolistsep]
        \item $e_{i-2},e_{i-1},e_i,e_{i+1}, e_{i+2}$ are all canonical in~$C$;
        \item $e_{i-1},e_i,e_{i+1}$ are all strictly canonical in~$C$.
    \end{itemize}
\end{lemma}
\begin{proof}
  If $C$ is fully canonical and no two anchors coincide, we can pick any edge of~$C$ as~$e_i$ and are done.
  Otherwise, we pick a sequence of canonical edges~$e_1,\dots,e_m$ in~$C$ such that~$v_1 = v_m$ and such that the spikes~$\Sp_2, \dots, \Sp_{m-1}$ are contained in the region bounded by~$\Lo \cup \Sp_1 \cup \Sp_m$.  (When $C$ is not fully canonical, then its sub-sequence~$e_2,\dots,e_{k-2}$ has this property by Theorem~\ref{thm:structure}.)
  
  Since spikes cannot intersect (except for touching at their anchors), the coinciding anchors form a bracket structure. Pick any innermost interval~$e_j, \dots, e_\ell$ such that~$v_j = v_\ell$. Then the anchors~$v_{j},\dots,v_{\ell-1}$ are all distinct. By Lemma~\ref{lem:strict}, $\ell - j \geq 5$. We can thus pick~$i = j + 2$ to satisfy the requirements of~the~lemma. \qed
\end{proof}

For each chordless odd cycle~$C$ of~$\IG$, we pick an edge~$e$ as in Lemma~\ref{lem:ground} and call it the \emph{ground edge} of~$C$.  Our proof of Theorem~\ref{thm:main} relies on the following key lemma:
\begin{lemma}
  \label{lem:key}
    Let $C$ and $C'$ be two chordless odd cycles with ground edges~$e$ and~$e'$.
    If~$e$ and~$e'$ cross, then $e$ is part of~$C'$ and $e'$ is part of~$C$.
\end{lemma}

\ifarxiv
The proof can be found in Appendix~\ref{sec:key}.
\fi
We can now prove our main theorem.
\begin{proof}[of Theorem~\ref{thm:main}]
Starting with $S = \emptyset$, we consider all chordless odd cycles~$C$ in~$\IG$ one by one.
If $S$ does not already contain an edge of~$C$, we add the ground edge of~$C$ to~$S$.

We claim that the resulting set~$S$ satisfies the conditions of the theorem. Indeed, by construction $S$ contains an edge of every chordless odd cycle. Since every odd cycle contains a chordless odd cycle as a subset, condition~(2) holds.  Assume now for a contradiction that there are two edges~$e, e' \in S$ such that~$e$ and~$e'$ cross.  Let $C$ and $C'$ be the chordless odd cycles that caused $e$ and~$e'$ to be added to~$S$, and assume w.l.o.g.~that $e$ was added to~$S$ before~$e'$.  By Lemma~\ref{lem:key}, $e$ is an edge of~$C'$---but that means that when considering~$C'$, no edge has been added to~$S$, a contradiction. \qed
\end{proof}

\section{Open problems}
\label{sec:openproblems}

We conclude with some open problems:
\begin{enumerate}
\item Do these results hold for fan-planar graphs in the ``classic'' definition?
\item Are there actually strongly fan-planar graphs that have thickness three?
\item Are there $2$-planar graphs of thickness three? 
\item What is the thickness of the optimal $2$-planar graphs? 
(These graphs have been fully characterized~\cite{LIPIcs:2017:7230}.)
It can be shown that these graph admit an edge decomposition into a $1$-planar graph and a bounded-degree planar graph~\cite{DBLP:journals/dm/BekosGDLMR19}. 
If the thickness of this graph is actually three, this would answer both previous questions in the affirmative as the optimal $2$-planar graphs are also fan-planar.
\item If there is no $2$-planar graph of thickness three, what is the smallest~$k$ such that there exists a $k$-planar graph of thickness three? 
Note that $K_9$ is $4$-planar and requires thickness three, hence $k \in \{2,3,4\}$. 
\item Are there strongly fan-planar graphs~$\G$ such that \emph{every} fan-planar drawing of~$\G$ requires three colors for the edges?  
In other words, are there strongly fan-planar graphs where odd cycles in the intersection graph of its drawing are unavoidable?
\end{enumerate}

\bibliographystyle{splncs04}
\bibliography{bibliography}

\ifarxiv
\newpage
\appendix
\section{Proof of Lemma~\ref{lem:e-cross-l}}
\label{sec:e-cross-l}

We repeat the statement of Lemma~\ref{lem:e-cross-l}:
\setcounter{lemma}{5}
\begin{lemma}
 Let~$C$ be a chordless cycle of length at least five, 
 and let~$e$ be an edge not part of~$C$ such that~$e$ intersects the loop~$\Lo$ of~$C$. 
 Then~$e$ starts in the anchor~$v_i$ of a canonical edge~$e_i$ of~$C$, 
 passes through the spike~$\Sp_i$, 
 crosses the base~$\base_i$, and either 
 (1)~ends in~$\Lo$; 
 or (2)~crosses the base~$\base_j$ of another canonical edge~$e_j$ of~$C$, 
 then passes through~$\Sp_j$ and terminates in~$v_j$; or 
 (3)~crosses the base of the non-canonical edge~$e_j$ of~$C$, never enters~$\Lo$ again, 
 and terminates in a vertex that is not a vertex of~$C$. 
 The second and third case can only happen if~$e_i$ and~$e_j$ share an endpoint. 
 In particular, (2) cannot happen when~$e_{i-1}$ and~$e_{i+1}$ are strictly canonical, 
 and (3) implies that~$i \in \{j-2, j+2\}$.
\end{lemma}
\setcounter{lemma}{8}
\begin{proof}
We first observe the following: 
if~$e$ intersects~$\Lo$ in a point~$p$ on a canonical edge~$e_i$, so that~$p \in \base_i$,
then by~\reffani the anchor~$v_i$ must be an endpoint of~$e$.  
By~\reffanii and~\refsimple, 
this means that the part of~$e$ between~$v_i$ and~$p$ lies in the spike~$\Sp_i$.

Assume first that~$e$ does not cross a non-canonical edge of~$C$. 
If~$e$ crosses~$\Lo$ once (in canonical edge~$e_i$), then by the above we are in case~(1):
$e$~starts in~$v_i$, passes through~$\Sp_i$,
crosses~$\base_i$, and ends inside~$\Lo$.
If~$e$ crosses~$\Lo$ more than once, 
then by the above this must be in two canonical edges~$e_i$ and~$e_j$.
Then we are in case~(2): 
$e$ starts in~$v_i$, passes through~$\Sp_i$, crosses~$\base_i$, 
passes through~$\Lo$, crosses base~$\base_j$, passes through~$\Sp_j$, and ends in~$v_j$.
Since~$e$ intersects~$e_i$ and~$e_j$, 
\reffani implies that~$e_i$ and~$e_j$ have a common endpoint.
If $e_{i-1}$ and $e_{i+1}$ are strictly canonical, 
then the endpoints of~$e_i$ are the anchors of~$e_{i-1}$ and~$e_{i+1}$
and both have degree two in~$\G_C$, so~$j \in \{i-2, i+2\}$. 
But then~$e$ is identical to either~$e_{i-1}$ or~$e_{i+1}$, a contradiction.  

Next, we consider the case where~$e$ crosses the non-canonical edge~$e_j$ of~$C$ in a point~$q$.
By~\reffani, then the anchor~$v_j$ must be an endpoint of~$e$,
and by~\reffanii, the part of~$e$ connecting~$v_j$ and~$q$ must lie in the region 
bounded by~$e_{j-1}$, $\base_j$, and~$e_{j+1}$.  
But that means that~$e$ passes through~$\Lo$ on this part, 
so $e$~must cross another, canonical edge~$e_i$ of~$C$.
By the observation above, this means that~$v_i = v_j$, 
so~$e$ starts in~$v_i$, passes through~$\Sp_i$, crosses~$\base_i$,
passes through~$\Lo$, and crosses~$\base_j$.
Since~$e$ cannot cross~$\base_j$ again and, by the observation above,
cannot cross another canonical edge, $e$~does not enter~$\Lo$ again.
Since~$e$ intersects~$e_i$ and~$e_j$, 
\reffani implies that~$e_i$ and~$e_j$ have a common endpoint.
The only edges that share an endpoint with the non-canonical edge~$e_j$ 
are the edges~$e_{j-2}$ and~$e_{j+2}$, so we have~$i \in \{j-2, j+2\}$.
Finally, the other endpoint of~$e$ cannot be a vertex of~$C$:
the only vertices of~$C$ that~$e$ could reach now are the endpoints of~$e_j$ 
(a contradiction to~\refsimple), 
the anchor~$v_j$ (but then~$e$ is a loop),
or the endpoints~$b_{j-1}$ and~$a_{j+1}$.
But since~$v_i = a_{j-1} = b_{j+1}$, then~$e$ is identical to either~$e_{j-1}$ or~$e_{j+1}$, a contradiction.  
\qed
\end{proof}

\newpage
\section{Proof of Lemma~\ref{lem:key}}
\label{sec:key}

The proof of Lemma~\ref{lem:key} relies on the following two lemmas.
\begin{lemma}
  \label{lem:not-interior}
  Let $C$ and $C'$ be two chordless odd cycles with ground edges~$e$ and~$e'$ and such that~$e$ and~$e'$ cross in a point~$z$.  Then~$z$ does not lie in the relative interior of the base~$\base$.
\end{lemma}
\begin{proof}
We denote the edges of~$C$ as~$e_1,\dots,e_k$ and the edges of~$C'$ as~$e'_1,\dots,e'_\ell$. By renumbering the edges, we can assume that~$e = e_3$ and~$e' = e'_3$.  
Assume for a contradiction that~$z$ lies in the interior of~$\base_3$. 
By Lemma~\ref{lem:e-cross-l}, that means that~$e'_3$ starts in~$v_3$, passes through~$\Sp_3$, crosses~$\base_3$, 
and terminates at a point~$u$ inside~$\Lo$ without crossing it again (by the choice of ground edge~$e_3$, 
cases (2) and~(3) of Lemma~\ref{lem:e-cross-l} cannot occur).  
Since $e_3$ crosses~$e'_3$, one endpoint of~$e_3$ must be~$v'_3$.  Breaking symmetry, we assume~$b_3 = v_4 = v'_3$.  
We can also assume that~$v_3$ is the source~$a'_3$ of~$e'_3$ (otherwise we can reverse the orientation of~$C'$).
Consider~$e'_4$, which starts in~$v'_3 = v_4$. 
Edges~$e'_3$ and~$e'_4$ intersect in the point~$x'_4$ that lies on~$e'_3$ and therefore either in~$\Lo$ or in~$\Sp_3$. 
In the former case, $e'_4$ must intersect~$e_4$ to enter~$\Lo$ by Lemma~\ref{lem:e-cross-l}.  
In the latter case, we observe that ~$e'_4$ cannot enter~$\Sp_3$ by crossing~$e_2$ by~\reffanii, hence in order to enter~$\Sp_3$ it must also cross~$e_4$.  
The crossing of~$e_4$ and $e'_4$ implies that the anchor~$v'_4$ is an endpoint of~$e_4$.  
Since $v_3 = v'_2 \neq v'_4$, we have~$e_4 = (v'_2, v'_4) = e'_3$---but that means that~$z = e_3 \cap e'_3 = e_3 \cap e_4 = x_4$, 
and $z$~does not lie in the relative interior of~$\base$, a contradiction. \qed
\end{proof}

\begin{lemma}
  \label{lem:in-bases}
  Let $C$ and $C'$ be two chordless odd cycles with ground edges~$e$ and~$e'$ and such that~$e$ and~$e'$ cross in a point~$z$.  Then~$z$ lies in base~$\base$.
\end{lemma}
Before we prove Lemma~\ref{lem:in-bases}, we show how it implies Lemma~\ref{lem:key}.
Again we repeat the statement of the lemma first:
\setcounter{lemma}{7}
\begin{lemma}
    Let $C$ and $C'$ be two chordless odd cycles with ground edges~$e$ and~$e'$.
    If~$e$ and~$e'$ cross, then $e$ is part of~$C'$ and $e'$ is part of~$C$.
\end{lemma}
\setcounter{lemma}{10}
\begin{proof}
We use the notation of Lemma~\ref{lem:not-interior}, and assume that~$e_3$ and~$e'_3$ cross in a
point~$z$.  
By Lemma~\ref{lem:in-bases}, $z \in \base = \base_3$. By Lemma~\ref{lem:not-interior}, 
$z$~is then an endpoint of~$\base_3$, which implies that~$e'_3 = e_2$ or~$e'_3 = e_4$, 
so~$e' \in C$.
By symmetry, $e \in C'$. \qed
\end{proof}

The proof of Lemma~\ref{lem:in-bases} requires yet another lemma.
Fig.~\ref{fig:remain-spike} illustrates its rather technical conditions. 
\begin{figure}[htb]
  \centerline{\includegraphics[scale=0.95,page=11]{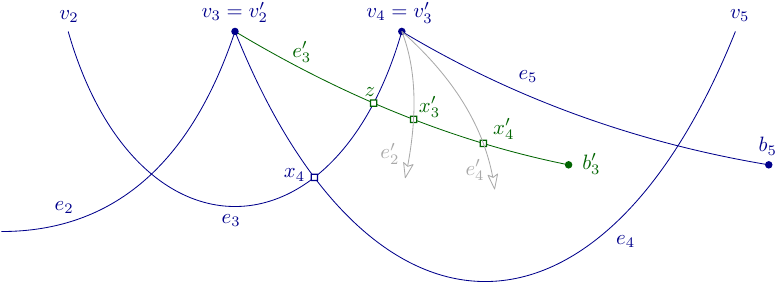}}
  \caption{Proof of Lemma~\ref{lem:remain-spike}.}
  \label{fig:remain-spike}
\end{figure}
\begin{lemma}
  \label{lem:remain-spike}
  Let $C$ and $C'$ be two chordless odd cycles with ground edges~$e_3$ and~$e'_3$. Suppose that $v_3 = v'_2$, $v_4 = v'_3$, that~$e'_3$ intersects~$x_4 v_4$ such that~$b'_3 \in \Sp_4$, and that~$x'_3 \in \Sp_4$.
  Then the entire loop~$\Lo'$ of~$C'$ lies in~$\Sp_4$.
\end{lemma}
\begin{proof}
We assume the contrary, so that~$\Lo'$ contains a point outside of~$\Sp_4$.
We follow the closed curve~$\Lo'$, starting from~$x'_4$, which lies strictly inside~$\Sp_4$.
There must be a first edge~$e'_m$ on~$C'$ such that~$x'_m$
lies still strictly inside~$\Sp_4$, but~$\base'_m$ intersects the boundary of~$\Sp_4$.
Recall that $\Sp_4$ is bounded by the edges~$e_3$, $e_4$ and $e_5$, so~$e'_m$ intersects
one of these edges, see Fig.~\ref{fig:remain-spike}.

Suppose first that $\base'_m$ crosses~$e_3$ in a point~$p$. 
Then~$e'_m$ is incident to~$v_3 = v'_2$.
Since~$e'_2$ is strictly canonical,~$v'_2$ has degree two in~$\G_{C'}$, and so~$e'_m = e'_1$.
Since~$e'_1$ does not intersect~$e_4$ or~$e'_3$ by~\refsimple, 
this means that~$p$ lies on~$e_3$ between~$x_4$ and~$z = e_3 \cap e'_3$,
and that~$x'_2$ lies in the triangle formed by~$e_4$, $e_3$, and~$e'_3$.
If $x'_2 = p$, then~$e'_2 = e_3$, but then~$x'_3 = z$ and~$\Lo'$ does not leave~$\Sp_4$.
If $x'_2$ lies strictly outside~$\Sp_4$,
then~$e'_2 \neq e_3$, and $x'_3$ lies strictly inside~$\Sp_4$.
But~$e'_2$ is incident to~$v'_3 = v_4$,
meets~$x'_3$ in the interior of~$\Sp_4$ and~$x'_2 \not\in \Sp_4$.
This is impossible because by~\refsimple~$e'_2$ does not intersect~$e_3$ or~$e_5$,
intersects~$e'_3$ only in~$x'_3$, and can intersect~$e_4$ at most once.

Suppose next that $\base'_m$ crosses~$e_4$. Then~$e'_m$ is incident to~$v_4 = v_3'$. 
But since $e_3'$ is strictly canonical, $v_3'$ has 
degree two in $\G_{C'}$---hence $e'_m = e'_4$.
It follows that~$e_4$ crosses~$e'_4$, so~$e_4$ is incident to~$v'_4$, 
and since~$v'_4 \neq v'_2$, 
this means that~$e_4 = (v'_2, v'_4) = e'_3$.  
But then~$b'_3 \not\in \Sp_4$, a contradiction.

Suppose last that $\base'_m$ crosses~$e_5$.  
Then~$v_m'$ is an endpoint of $e_5$, that is~$v'_m = a_5 = v_4 = v'_3$ or~$v'_m = b_5$. 
But~$e'_3$ is strictly canonical, so~$v'_3$ is the anchor of no other edge of~$C'$ than~$e'_3$,
so~$v'_m = v'_3$ implies~$e'_m = e'_3$, a contradiction.
Hence, $v_m' = b_5$. By construction~$x_m'$ is still contained strictly inside~$\Sp_4$, 
so~$e_{m-1}'$, which is incident to~$v'_m = b_5$ and contains~$x'_m$, 
has to intersect the boundary of $\Sp_4$. 
By~\refsimple $e_{m-1}'$ cannot intersect~$e_5$.
If~$e'_{m-1}$ crosses~$e_4$,
then~$e'_{m-1}$ is also incident to~$v_4 = v_3' = a_5$, that is~$e'_{m-1} = e_5$.
But~$e_5$ contains no point inside~$\Sp_4$, a contradiction.
Finally, by~\refminimal, $e'_{m-1}$ cannot cross~$e'_2$ or~$e'_3$, so even
if~$e'_{m-1}$ crosses~$e_3$, it can only reach the triangular region formed by
$e'_3$, $e'_2$, and~$e_3$, and~$\Lo'$ does not enter this region 
since~$\Lo'$ has no point on~$a'_3 x'_3 = v_3 x'_3$ or~$x'_3 b'_2 = x'_3 v_4$.
This final contradiction concludes the proof. \qed
\end{proof}

\begin{figure}[t]
  \centerline{\includegraphics[scale=0.95,page=9]{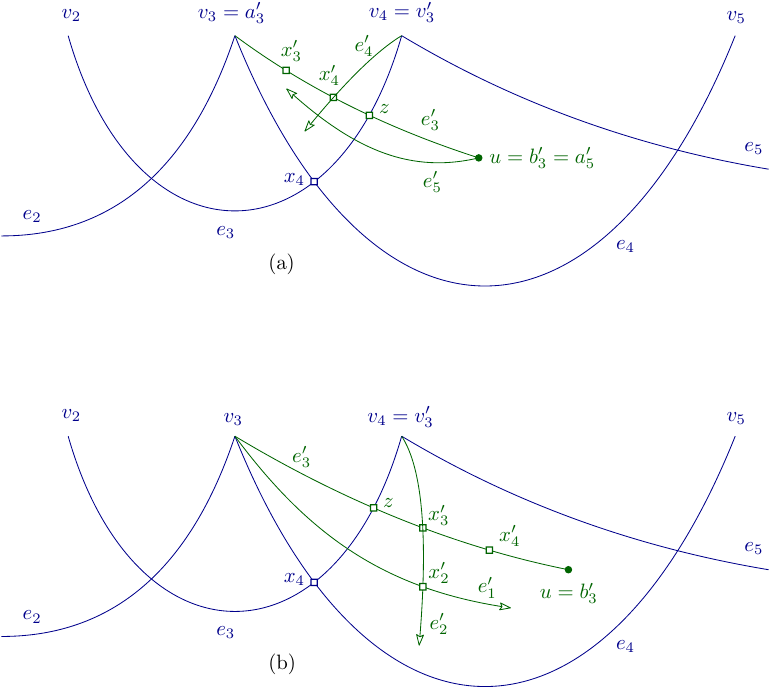}}
  \caption{Proof of Lemma~\ref{lem:in-bases}.}
  \label{fig:not-exterior}
\end{figure}
\begin{proof}[of Lemma~\ref{lem:in-bases}]
We use the notation of Lemma~\ref{lem:not-interior} and assume that $z \not\in \base$. 
Breaking symmetry, we can assume that~$z \in x_4 b_3$.  
Since~$e'_3$ crosses~$e_3$, $e'_3$ has an endpoint in~$v_3$.  
By Lemma~\ref{lem:vertex-is-in-spike}, this means that the other endpoint~$u$ of~$e'_3$ lies in~$\Sp_4$,
see Fig.~\ref{fig:not-exterior}(a).
Since $e_3$ crosses~$e'_3$, $v'_3$~is an endpoint of~$e_3 = (v_2, v_4)$.

If~$v'_3 = v_2$, then consider~$e'_2$ and~$e'_4$, which are both incident to~$v'_3 = v_2$.
By Lemma~\ref{lem:vertex-in-next-spike}, they cannot intersect the interior of~$\Sp_4$, 
so they can cross~$e'_3$ only on the segment~$v_3 z$.  
Without violating~\reffanii or~\reffaniii, 
this means they both need to cross~$e_4$.
But then both~$e'_2$ and~$e'_4$ must be incident to~$v_4$,
so then~$e'_2 = e'_4$, a contradiction.

It follows that~$v'_3 = v_4$. 
We can orient~$C'$ such that~$v_3$ is the source~$a'_3$ of~$e'_3$, 
and consider the order of~$x'_3$, $x'_4$,
and~$z$ on~$e'_3$. By Lemma~\ref{lem:not-interior}, 
either $z$ occurs after~$x'_3$ and~$x'_4$ (possibly with~$z = x'_4$), 
or before~$x'_3$ and~$x'_4$ (possibly with~$z = x'_3$).

If $x'_3$ and $x'_4$ are encountered before~$z$ (the situation shown in Fig.~\ref{fig:not-exterior}(a)), 
then consider the edge~$e'_5$. 
It cannot cross~$e'_3$~by~\refsimple,
and cannot cross~$e_4$, because then~$b'_5 = v_4 = v'_3 = a'_4$, a contradiction to~\refsimple 
since~$e'_4$ and~$e'_5$ cross.
Since~$e'_5$ crosses~$e'_4$ on the segment~$x'_4b'_4$, $e'_5$ must cross~$e_3$. 
Therefore~$e'_5$ is incident to~$v_3$---but then~$e'_5 = e'_3$, a contradiction. 

It follows that $x'_3$ and $x'_4$ are encountered after~$z$,
see Fig.~\ref{fig:not-exterior}(b).  We are now in the setting of Lemma~\ref{lem:remain-spike},
so the entire closed path~$\Lo'$ lies in~$\Sp_4$.  In particular,~$x'_2 \in \Sp_4$. 
But by~\refsimple, $e'_1$ cannot cross~$e_4$, 
and it cannot cross~$e_5$ by~\reffanii.  
It follows that~$e'_1$, which is incident to~$a'_3 = v_3$, must cross~$e_3$.
This then implies that $v'_1$ is an endpoint of~$e_3$.  
Since~$v_4 = v'_3 \neq v'_1$, this means that~$e_3 = (v'_1, v'_3) = e'_2$, 
and so $z = x'_3$ and~$x'_2$ lies on the segment of~$e_3$ between~$x_4$ and~$x'_3$,
\begin{figure}[htb]
  \centerline{\includegraphics[scale=0.95,page=7]{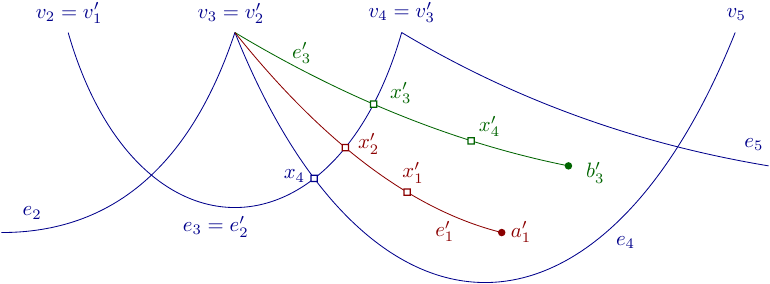}}
  \caption{Proof of Lemma~\ref{lem:in-bases}, having established~$e'_2 = e_3$.}
  \label{fig:in-bases}
\end{figure}
see Fig.~\ref{fig:in-bases}.
We now consider the edges~$e'_1$ and~$e'_0$.
Since~$b'_1 = v_3$, Lemma~\ref{lem:vertex-is-in-spike} 
implies that $a_1'$ is contained inside~$\Sp_4$,
and therefore~$x'_1$ is also strictly inside~$\Sp_4$.
Since~$e'_1$ is canonical, $e_0'$ is incident to $v_1' = a'_2 = v_2$.
But by Lemma~\ref{lem:vertex-in-next-spike}, then~$e'_0$ cannot intersect the interior
of~$\Sp_4$, so it cannot contain~$x'_1$, a contradiction. \qed
\end{proof}
\fi
\end{document}